\documentclass{article}%
\usepackage{amsfonts}
\usepackage{amsmath}
\usepackage{amsthm}
\usepackage{amssymb}
\usepackage{color}
\def\Q{{\mathbb Q}}

\def\th{{\theta}}

\def\z{{\zeta}}

\def\proofend{\hspace*{1mm} \hfill{$\Box$}}

\newcommand\jac[2]{\left(\frac{#1}{#2}\right)}
\newcommand{\froot}[1]{\frac{1-\z^{#1}}{1+\z^{#1}}}

\newtheorem{theorem}{Theorem}
\newtheorem{property}[theorem]{Property}

\begin{document}

\title{Elementary Trigonometric Sums \\ related to Quadratic Residues}
\author{A.~Laradji\thanks{Department of Mathematics \& Statistics, KFUPM, Dhahran,
Saudi Arabia, e-mail: alaradji@kfupm.edu.sa} 
\and
M.~Mignotte\thanks{Universit\'{e} Louis Pasteur,U.~F.~R. de Math\'{e}matiques, Strasbourg, France,
e-mail: maurice@math.u-strasbg.fr}
\and
N.~Tzanakis\thanks{Department of Mathematics, University of Crete, Iraklion,
Greece, e-mail: tzanakis@math.uoc.gr, http://www.math.uoc.gr/\~{}tzanakis}
}
\date{} 

\maketitle

Let $p$ be an odd prime. We will study the following sums
\begin{equation} \label{def T}
 T(p)=\sqrt{p}\sum_{n=1}^{(p-1)/2}\tan\frac{\pi n^2}{p}
\end{equation}
and
\begin{equation} \label{def C}
 C(p)=\sqrt{p}\sum_{n=1}^{(p-1)/2}\cot\frac{\pi n^2}{p}\,.
\end{equation}
Surprisingly, we came across these sums as we were working on a certain diophantine equation. Being 
non specialists in the relevant area, we were impressed by the nice properties that these sums have and their 
elegant consequences. It is this feeling of elegance that we would like to share with our readers.
As pointed out to us by Juan Carlos Peral Alonso, to whom we are grateful, these sums are closely related 
to the {\em class-number formula} due to Dirichlet (see (\ref{Lebesgue formula})), sometimes called 
``Lebesgue's formula'' --see \cite{Dick}, p.~179--, which ``explains'', in a sense, their nice properties.
For those readers who are not already acquainted with the notion of {\em class-number}, a brief remark 
has its place. Let $D$ be a negative integer which is a {\em fundamental discriminant}, i.e. either
$D\equiv 1\pmod 4$ and $D$ is squarefree, or $D\equiv 0\pmod 4$ and $D/4$ is squarefree $\equiv 2,3\pmod 4$.
In particular, if $p$ is a prime $\equiv 3\bmod 4$ (we will deal with such primes in this paper), 
$-p$ is a fundamental discriminant. The {\em class-number} $h(D)$ has a double interpretation, as the
number of reduced binary quadratic forms of discriminant $D$, and as the number of classes of fractional
ideals of the quadratic number field $\Q(\sqrt{D})$. The reader may very well profit by reading, for example,
sections 4.9.1, 5.1, 5.2 and 5.3.1 of H.~Cohen's book \cite{Cohen}, written in a very concrete way; see,
especially, the conclusion following Lemma 5.3.4 therein.

All the results presented in this paper, possibly with the exception of Properties  \ref{property 1},
\ref{property 3} and \ref{property h}(\ref{h and odd and even squares}), are scattered in the literature, mainly 
(but not exclusively) in articles about the class-number of binary quadratic forms; 
see, for example, \cite{Dick} and \cite{Lerch}. 
Therefore, our purpose is not to present new results; rather having expository-pedagogic aim, our paper 
offers a bouquet of classical results which are presented in a very smooth, as we believe, manner, 
practically using only Elementary Mathematics, or appealing to short and easily readable elementary
papers, like \cite{BaCho}, \cite{BerChow}, \cite{Chung},\cite{Mo}, \cite{Whi}.

Since $T(p)$ and $C(p)$ are very closely related to each other (see (\ref{relation of T and M}) and (\ref{relation of C and M})), 
we will mainly focus on $T(p)$. We also note that our $T(p)$ is equal to H.L.~Montgomery's
$-T(1,\chi)$ as defined in \cite{Montg}, where $\chi$ is the non-trivial quadratic character.

As we will see immediately below, if $p\equiv 1\pmod 4$, then $T(p)=0=C(p)$, therefore, concerning
the sums $T(p)$ and $C(p)$, only the case $p\equiv 3\pmod 4$ is of interest. For this case,
we prove a number of elegant number-theoretical properties of $T(p)$.
Some of them have the flavour of the well-known property of the primes 
$p\equiv 3\pmod 4$, asserting that, in the range $1$ to $(p-1)/2$ there are more quadratic residues $\bmod\,p$ 
than non-quadratic residues (see, for example, \cite{Chung}, \cite{Mo}, \cite{Whi}). 
Further, Properties \ref{property 1} and \ref{property 2} below give a simple rule comparing the numbers of even 
and odd quadratic residues in $\{1,2,\ldots,p-1\}$ and Property \ref{property h} gives an extremely simple
rule for expressing $h(-p)$, the class-number of the quadratic field $\Q(\sqrt{-p})$. 

First, a few remarks have their place.
If $Q$ is any complete set of quadratic residues $\bmod\,p$, we can write
\[
 T(p)=\sqrt{p}\sum_{j \in Q}\tan\frac{\pi j}{p}\,.
\]
If $p\equiv 1\pmod 4$, then $-Q\equiv Q\pmod p$, from which we immediately conclude that $T(p)=0$
and, similarly, $C(p)=0$. Therefore we make the following assumption:
\begin{quote}
 Throughout this paper, $p$ will always denote a prime $\equiv 3\pmod 4$.
\end{quote}
We denote by $\z$ a primitive $p$-root of unity and we put $i=\sqrt{-1}$. 
Also, by $\sqrt{p}$ we mean the positive square root of $p$.
It is easy to see that
\begin{equation} \label{T second form}
 T(p)=i\sqrt{p}\sum_{j\in Q}\frac{1-\z^j}{1+\z^j} \,,
\end{equation}
therefore, we have
\begin{align}
 T(p) = & i\sqrt{p}\sum_{j=1}^{(p-1)/2}\froot{j^2} = 
i\sqrt{p}\sum_{j=1}^{(p-1)/2}\left(\frac{2}{1+\z^{j^2}}-1\right) =
i\sqrt{p}\sum_{j=1}^{(p-1)/2}\left(\frac{1+(\z^{j^2})^p}{1+\z^{j^2}}-1\right) 
                \nonumber \\
   =  & i\sqrt{p}\sum_{j=1}^{(p-1)/2}\left(1-\z^{j^2}+(\z^{j^2})^2-\cdots +(\z^{j^2})^{p-1}-1\right) 
               \nonumber \\
 = & i\sqrt{p}\sum_{j=1}^{(p-1)/2}\left(-\z^{j^2}+(\z^{j^2})^2-\cdots +(\z^{j^2})^{p-1}\right)  
                   \nonumber \\
 = & \frac{i\sqrt{p}}{2}\,\sum_{j=1}^{p-1}\left(-\z^{j^2}+(\z^{j^2})^2-\cdots +(\z^{j^2})^{p-1}\right) 
                 \nonumber \\
 = & \frac{i\sqrt{p}}{2}\,\sum_{k=1}^{p-1}(-1)^k\sum_{j=1}^{p-1}\z^{j^2k}= 
                              \frac{i\sqrt{p}}{2}\,\sum_{k=1}^{p-1}(-1)^k\sum_{j=0}^{p-1}\z^{j^2k}\,.
              \label{eq T as Gauss}
\end{align}
For every $k=1,\ldots, p-1$, $\displaystyle{\sum_{j=0}^{p-1}\z^{j^2k}}$ is the well-known
Gaussian sum, denoted by $S(k,p)$, which is equal to $i\jac{k}{p}\sqrt{p}$, where $\jac{\cdot}{p}$ is the 
Legendre symbol. This is a straightforward consequence of the following more general well-known result:
\\
{\em Let $m$ be an odd positive number and let $n$ be an integer relatively prime to $m$. \\
Put $\displaystyle{ S(k,m)=\sum_{j=0}^{m-1}e^{2\pi ij^2k/m}}$. Then,
\[
 S(k,m)=\begin{cases}
         \jac{k}{m}\sqrt{m} & \mbox{if $m\equiv 1\pmod 4$} \\[1mm]
         i\jac{k}{m}\sqrt{m} & \mbox{if $m\equiv 3\pmod 4$}
        \end{cases} 
\]
}
See, e.g. Theorem 5.6 in Chapter 7 of \cite{Hua}. When $m$ is a prime $p\equiv 3\pmod 4$,
we can more directly prove that 
\begin{equation} \label{eq Gauss Sum}
S(k,p)=i\jac{k}{p}\sqrt{p}\,,
\end{equation}
without appealing to the above result, by turning to a short paper of Bamba and Chowla \cite{BaCho}. 
In that paper, an interesting brief and elementary proof of the relation 
$ (1-i)(1+i^m)S(1,m)=2\sqrt{m}$, where $m$ is positive odd integer, is given. 
Consequently, if $p$ is a prime $\equiv 3\pmod 4$, then 
$S(1,p)=i\sqrt{p}$. By the definition of $S(k,m)$ it is clear that, if $k$ is a quadratic residue $\bmod\,p$, 
then $S(k,p)=S(1,p)$; and if $k$ is a quadratic non-residue, then 
\[
 S(k,p)=S(-1,p)=\overline{S(1,p)}= \overline{i\sqrt{p}}=-i\sqrt{p}\,,
\]
as claimed.
\\
Now, going back to (\ref{eq T as Gauss}) and using (\ref{eq Gauss Sum}), we obtain the following expression 
for $T(p)$:
\begin{equation} \label{first expression T}
 T(p)=\frac{p}{2}\,\sum_{k=1}^{p-1}(-1)^{k+1}\jac{k}{p} \,.
\end{equation}
Let $a_1,a_2,\ldots,a_{\mu}$ and $b_1,b_2,\ldots,b_{\nu}$ be, respectively, the even and odd quadratic
residues $\bmod\,p$ in the set $P=\{1,2,\ldots,p-1\}$. Clearly, $\mu+\nu=(p-1)/2$ and the set of the quadratic
non-residues $\bmod\,p$ in $P$ is $\{p-a_1,\ldots p-a_{\mu},p-b_1,\ldots,p-b_{\nu}\}$. Note that a summand
$(-1)^{k+1}\jac{k}{p}$ in the right-hand side of (\ref{first expression T}) is positive iff
$k\in\{p-b_1,\ldots,p-b_{\nu}, b_1,\ldots,b_{\nu}\}$, i.e. $2\nu$ summands are positive and, analogously, $2\mu$
summands are negative. Then, $T(p)=p(\nu-\mu)$, where we observe that $\nu-\mu$ is an odd number, since
$\nu+\mu=(p-1)/2$. Thus, we have the following:
\begin{property} \label{property 1}
Let $p$ be a prime $\equiv 3\pmod 4$ and let $q_o(p)$ and $q_e(p)$ be, respectively, the number of odd 
and even quadratic residues $\bmod\,p$ in the set $\{1,2,\ldots,p-1\}$. Then
\begin{equation} \label{second expression T}
 T(p)=p(q_o(p)-q_e(p)) \,.
\end{equation}
In particular, $T(p)$ is an odd integer divisible by $p$ and by no higher power of $p$.
\end{property}
Next, we rewrite the definition (\ref{def T}) of $T(p)$ as follows, 
\begin{equation} \label{T third form}
 T(p)=\frac{\sqrt{p}}{2}\,\sum_{n=1}^{p-1}\tan\frac{n^2\pi}{p}\,.
\end{equation}
We have the following inequality of A.L.~Whiteman (Theorem 2 of \cite{Whi}):
\begin{equation} \label{Whi ineq}
 \sum_{n=1}^{p-1}\cot\frac{n^2\pi}{p} > 0\,.
\end{equation}
In view of the identity $\tan\th=\cot\th-2\cot 2\th$, the relation (\ref{T third form}) becomes
\begin{equation} \label{sum tan}
 \frac{2}{\sqrt{p}}T(p) =\sum_{n=1}^{p-1}\cot\frac{n^2\pi}{p}-2\sum_{n=1}^{p-1}\cot\frac{2n^2\pi}{p}\,.
\end{equation}
Let $n=1,2,\ldots,p-1$.  
If $p\equiv 7\pmod 8$, the sets $\{2n^2\}$ and $\{n^2\}$ are identical $\bmod\,p$, hence the right-hand 
side of (\ref{sum tan}) is equal to $\displaystyle{-\sum_{n=1}^{p-1}\cot\frac{n^2\pi}{p}}$ and, by Whiteman's 
inequality (\ref{Whi ineq}), it is strictly negative.
If $p\equiv 3\pmod{8}$, the sets $\{-2n^2\}$ and $\{n^2\}$ are identical $\bmod\,p$, therefore, the 
right-hand side of (\ref{sum tan}) is equal to $\displaystyle{3\sum_{n=1}^{p-1}\cot\frac{n^2\pi}{p}}$, 
hence, by (\ref{Whi ineq}), it is strictly positive. 
Thus, in combination also with Property \ref{property 1}, we obtain the following:
\begin{property} \label{property 2}
 $T(p)>0$ if $p\equiv 3\pmod 8$ and $T(p)<0$ if $p\equiv 7\pmod 8$. Also, in the set $\{1,2,\ldots,p-1\}$,
the odd quadratic residues $\bmod\,p$ are more than the even ones when $p\equiv 3\pmod 8$; 
the reverse situation is true when $p\equiv 7\pmod 8$.
\end{property}
Now consider the sum
\[
 M(p)=-\sum_{k=1}^{p-1}\jac{k}{p}k\,.
\]
Dirichlet \cite{Dir1} proved that, for $p\equiv 3\pmod 4$, $M(p)>0$, i.e. among the numbers
$1,2,\ldots, p-1$, the sum of the quadratic non-residues is greater than the sum of the quadratic residues.
In \cite{Ber}, B.C.~Berndt proves that 
\begin{equation} \label{another expression of M}
 M(p) =\frac{\sqrt{p}}{2}\sum_{k=1}^{p-1}\jac{k}{p}\cot\frac{k\pi}{p}
\end{equation}
and, based on (\ref{another expression of M}), he gives (Theorem 3.1 in \cite{Ber}) another 
proof of Dirichlet's inequality 
\begin{equation} \label{M>0}
 M(p) > 0 \quad\mbox{for $p\equiv 3\pmod 4$.}
\end{equation}
Using (\ref{another expression of M}), it is an easy exercise to check that
\begin{equation} \label{relation of T and M}
 T(p)=\begin{cases}
       -M(p) & \mbox{if $p\equiv 7\pmod 8$} \\
       3M(p) &  \mbox{if $p\equiv 3\pmod 8$} 
      \end{cases}
\end{equation}
This, combined with Property \ref{property 2}, gives now another proof of (\ref{M>0}) which is
simpler than that of Theorem 3.1 in \cite{Ber}.

{\bf An upper bound for} $T(p)$. From (\ref{second expression T}) we trivially obtain  
$|T(p)|<p(p-1)/2$.
However, we can obtain a much better upper bound as follows.
\\
Let $Q\subseteq \{1,2,\ldots,p-1\}$ be a complete set of quadratic residues $\bmod\,p$. We have
\[
 |T(p)|\leq \sqrt{p}\sum_{j\in Q}\left|\tan\frac{\pi j}{p}\right| =
          \sqrt{p}\sum_{j\in Q}\frac{1}{\left|\tan\frac{\pi(p-2j)}{2p}\right|}\,.
\]
Since $\left|\frac{\pi(p-2j)}{2p}\right|<\frac{\pi}{2}$, it follows that 
$\left|\tan\frac{\pi(p-2j)}{2p}\right|>\frac{\pi}{2p}|p-2j|$, hence,
\[
 |T(p)| <\frac{2p\sqrt{p}}{\pi}\sum_{j \in Q} \frac{1}{|p-2j|}\,.
\]
Note that, as $j$ runs through the set $Q$, the numbers $|p-2j|$ are distinct $\bmod\,p$, for, if
$|p-2j_1|\equiv |p-2j_2| \pmod p$ with $j_1,j_2\in Q$ and $j_1\neq j_2$, then, necessarily, $j_2\equiv -j_1\pmod p$,
which implies that $-1$ is a quadratic residue $\bmod\,p$, a contradiction. Therefore, the set 
$\{|p-2j|: j\in Q\}$ is a subset of $\{1,\ldots,p-1\}$ with cardinality $(p-1)/2$, consisting of odd numbers,
i.~e. it coincides with $\{1,3,\ldots,p-2\}$. Therefore, 
\[
 \sum_{j \in Q} \frac{1}{|p-2j|} \leq \sum_{k=1}^{(p-1)/2}\frac{1}{2k-1} < 1+\frac{1}{2}\log(p-2)\,,
\]
from which we obtain the following:
\begin{property} \label{property 3}
For any prime $p \equiv 3\pmod 4$ we have
\begin{equation} \label{ub for T}
|T(p)| < \frac{2p\sqrt{p}}{\pi}\left(1+\frac{1}{2}\log(p-2)\right).
\end{equation}
\end{property}
Now we go on to the study of $C(p)$. We use the following alternative
expression for $C(p)$ (cf. (\ref{T second form})):
\begin{equation} \label{C second form}
 C(p)=-i\sqrt{p}\sum_{j\in Q}\frac{1+\z^j}{1-\z^j} \,,
\end{equation}
where $Q$ is a complete set of quadratic residues $\bmod\,p$.
It is straightforward to check that $C(3)=1$, therefore we assume that $p>3$.
By Whiteman's inequality (\ref{Whi ineq}), we have $C(p)>0$.
\\
Just before obtaining Property \ref{property 2}, we actually proved that
$T(p)=-C(p)$ if $p\equiv 7\pmod 8$ and $T(p)=3C(p)$ if $p\equiv 3\pmod 8$. Therefore,
in view of (\ref{relation of T and M}),
\begin{equation} \label{relation of C and M}
C(p)=M(p).
\end{equation}
Relations (\ref{relation of T and M}) and (\ref{relation of C and M}) combined with Property \ref{property 1}, imply the following:
\begin{property} \label{property 4}
$C(p)$ is equal to the excess of the sum of quadratic non-residues over the sum of quadratic residues $\bmod\,p$.
$C(p)$ is an odd positive integer, divisible by $p$ and by no higher power of $p$. 
If $p\equiv 3\pmod 8$, then $T(p)$ is a multiple of $3$, hence, $q_o(p)-q_e(p)$ is a positive multiple of $3$.
\end{property}

The following elegant property relates $T(p)$ with the class-number of the quadratic number field 
$\Q(\sqrt{-p})$.
\begin{property} \label{property h}
Let $p$ be a prime number $\equiv 3\pmod 4$ and let $h(-p)$ be the class-number of 
the quadratic number field $\Q(\sqrt{-p})$. Then,
\begin{equation} \label{h and T}
 T(p)= \begin{cases}
        -ph(-p) & \mbox{if $p\equiv 7\pmod 8$} \\
        3ph(-p) & \mbox{if $p\equiv 3\pmod 8$}
       \end{cases}
\end{equation}
\begin{equation} \label{h and C}
 M(p)=C(p)= ph(-p)
\end{equation}
\begin{equation} \label{h and odd and even squares}
 h(-p) = \begin{cases}
          q_e(p)-q_o(p) & \mbox{if $p\equiv 7\pmod 8$} \\
          \frac{1}{3}(q_o(p)-q_e(p))& \mbox{if $p\equiv 3\pmod 8$}
         \end{cases}
\end{equation}
\end{property}
\noindent
{\bf Proof}. We have 
\begin{equation} \label{Lebesgue formula}
 h(-p) = \frac{1}{2\sqrt{p}}\sum_{k=1}^{p-1}\jac{k}{p}\cot\frac{k\pi}{p}\,.
\end{equation}
This is a consequence of the more general formula, referred to as ``Lebesgue's formula'' in Dickson's 
``History'' \cite{Dick}, p.~179, due to Dirichlet \cite{Dir2}. For a recent proof of that formula we refer the reader
to the Corollary 2.3 of \cite{BerZah}.
\\
A complete set of quadratic non-residues $\bmod p$ is $\{-k^2:\,k=1,\ldots,(p-1)/2\}$. Therefore, 
(\ref{Lebesgue formula}) becomes
\begin{align*}
 h(-p) & =\frac{1}{2\sqrt{p}}
\left(\sum_{k=1}^{(p-1)/2}\cot\frac{k^2\pi}{p}-\sum_{k=1}^{(p-1)/2}\cot\frac{-k^2\pi}{p}\right) \\
  & = \frac{1}{\sqrt{p}}\sum_{k=1}^{(p-1)/2}\cot\frac{k^2\pi}{p}=\frac{1}{p}C(p)=
                                   \mbox{(by (\ref{relation of C and M}))\:} \frac{1}{p}M(p)\,,
\end{align*}
which proves (\ref{h and C}), and now (\ref{h and T}) and (\ref{h and odd and even squares}) are straightforward 
consequences of (\ref{relation of T and M}) and Property \ref{property 1}, respectively.
\proofend
\\
The relation (\ref{h and T}) is a special case of Corollary 5.2 in \cite{BerZah} which goes back to
V.A.~Lebesgue \cite{Leb}.
The relation (\ref{h and C}) is due to Dirichlet \cite{Dir2}; see also Corollary 3.6 of \cite{BerZah}.
\\
Note that, since $q_e(p)+q_o(p)=(p-1)/2$, which is odd, Property \ref{property h} implies the following:
\begin{center}
{\em For a prime $p\equiv 3\pmod 4$, $h(-p)$ is odd.}
\end{center}
This is Corollary 3.6 of \cite{Ber}.
 
{\bf Further expressions for $T(p)$ and consequences}. Since $\jac{p-k}{p}=-\jac{k}{p}$, 
we have from (\ref{first expression T}),
\begin{equation} \label{first expression up to p/2}
  T(p)=p\sum_{k=1}^{(p-1)/2}(-1)^{k+1}\jac{k}{p} \,.
\end{equation}
We have
\begin{align}
\sum_{k=1}^{p-1}(-1)^{k+1}\jac{k}{p} & = 
\jac{1}{p}+\jac{3}{p}+\cdots +\jac{p-2}{p}-\jac{2}{p}-\jac{4}{p}-\cdots -\jac{p-1}{p} 
                              \label{odd numer minus even numer} \\
 & =\jac{1}{p}+\jac{3}{p}+\cdots +\jac{p-2}{p}
+\jac{p-2}{p}+\jac{p-4}{p}+\cdots +\jac{1}{p} \nonumber \\
& = 2\left(\jac{1}{p}+\jac{3}{p}+\cdots +\jac{p-2}{p}\right) = 2A\,, \label{2A}
\end{align}
where $A=\jac{1}{p}+\jac{3}{p}+\cdots +\jac{p-2}{p}$. Therefore, by (\ref{2A}) and 
(\ref{odd numer minus even numer}),
\[
2A= A-\jac{2}{p}\left(\jac{1}{p}+\jac{2}{p}+\cdots \jac{(p-1)/2}{p}\right)\,,
\]
implying
\[
 A=-\jac{2}{p}\left(\jac{1}{p}+\jac{2}{p}+\cdots \jac{(p-1)/2}{p}\right)\,.
\]
Collecting together the expressions for $T(p)$ in (\ref{first expression T}), (\ref{first expression up to p/2}),
(\ref{2A}) and the expression for $A$ just above, we have:
\begin{align}
 T(p) & = \frac{p}{2}\sum_{k=1}^{p-1}(-1)^{k+1}\jac{k}{p}   \label{T1}    \\
      & = p\sum_{k=1}^{(p-1)/2}(-1)^{k+1}\jac{k}{p}      \label{T2}  \\
      & = p\left(\jac{1}{p}+\jac{3}{p}+\cdots +\jac{p-2}{p}\right) \label{T3}  \\
      & = -p\jac{2}{p}\left(\jac{1}{p}+\jac{2}{p}+\cdots \jac{(p-1)/2}{p}\right)\,. \label{T4} 
\end{align}
Note that the sum appearing in (\ref{T4}) is a sum of the values of a primitive character $\bmod\,p$, 
therefore, by the P\'{o}lya inequality\footnote{Or P\'{o}lya-Vinogradov inequality.} 
(see Theorem 8.21 in \cite{Apo} or inequality (2) in Chapter 23 of \cite{Dav}) this sum is $< p^{1/2}\log p$. 
This gives the upper bound $|T(p)|<p^{3/2}\log p$, which is slightly worse than the upper bound 
in Property \ref{property 3}.

The expression (\ref{T4}) for $T(p)$, in combination with a well-known result of Dirichlet saying that,
among the numbers $1,2,\ldots, (p-1)/2$ there are more quadratic residues than non-quadratic residues $\bmod\,p$
(see e.g. \cite{Chung}, \cite{Whi}, \cite{Mo}, or exercises 14 through 17, Chapter 16 of \cite{IrRo})
furnishes another proof of Property \ref{property 2}.
\\
Next, equating the right-hand sides of (\ref{T2}) and (\ref{T4}) and separating the Legendre symbols with
even ``numerators'' from those with odd ones, we find 
\begin{align} 
&  \left(1+\jac{2}{p}\right)\left(\jac{1}{p}+\jac{3}{p}+\cdots +\jac{(p-1)/2}{p}\right) \nonumber \\
 & = \left(1-\jac{2}{p}\right)\left(\jac{2}{p}+\jac{4}{p}+\cdots +\jac{(p-3)/2}{p}\right)\,.
             \label{even and odd numerators}
\end{align}
If $p\equiv 3\pmod 8$, then (\ref{even and odd numerators}) implies that
$\jac{2}{p}+\jac{4}{p}+\cdots +\jac{(p-3)/2}{p}=0$, hence, $\jac{1}{p}+\jac{2}{p}+\cdots +\jac{(p-3)/4}{p}=0$,
which says that, among the numbers $1,2,\ldots,(p-3)/4$ there are as many quadratic residues as quadratic
non-residues $\bmod\,p$; and since $T(p)>0$, (\ref{T4}) implies that there are more quadratic residues than
quadratic non-residues among the numbers $(p+1)/4,\ldots,(p-1)/2$. 
\\
If $p\equiv 7\pmod 8$, then (\ref{even and odd numerators}) implies that 
$\jac{1}{p}+\jac{3}{p}+\cdots +\jac{(p-1)/2}{p}=0$, hence, 
$\jac{p-1}{p}+\jac{p-3}{p}+\cdots +\jac{(p+1)/2}{p}=0$ and, consequently,
$\jac{(p+1)/4}{p}+\cdots +\jac{(p-3)/2}{p}+\jac{(p-1)/2}{p}=0$. This shows that there are as many 
quadratic residues as quadratic non-residues $\bmod\,p$ among the numbers $(p+1)/4,\ldots,(p-1)/2$;
and since $T(p)<0$, (\ref{T4}) implies that there are more quadratic residues than
quadratic non-residues among the numbers $1,\ldots,(p-3)/4$. 
\begin{property} \label{property 5}
If $p\equiv 3\pmod 8$ then, among the numbers $1,2,\ldots,(p-3)/4$, there are as many quadratic residues 
as quadratic non-residues $\bmod\,p$ and among the numbers \\$(p+1)/4,\ldots,(p-1)/2$ the quadratic
residues are more than the quadratic non-residues.
\\
If $p\equiv 7\pmod 8$ then, among the numbers $1,2,\ldots,(p-3)/4$, there are more quadratic non-residues 
than quadratic residues $\bmod\,p$ and among the numbers $(p+1)/4,\ldots,(p-1)/2$ the quadratic
residues are as many as the quadratic non-residues.
\\
In other words,
\begin{eqnarray} 
\sum_{k=1}^{(p-3)/4}\jac{k}{p} & 
         \begin{cases}
              >0 & \mbox{if $p\equiv 7\pmod 8$} \\
             = 0 & \mbox{if $p\equiv 3\pmod 8$}
         \end{cases}                        \label{Berndt-Chowla 1}
\\
\sum_{k=(p+1)/4}^{(p-1)/2}\jac{k}{p} & 
         \begin{cases}
              =0 & \mbox{if $p\equiv 7\pmod 8$} \\
              >0 & \mbox{if $p\equiv 3\pmod 8$}
         \end{cases}                            \label{Berndt-Chowla 2}
\end{eqnarray}
\end{property}
The relations (\ref{Berndt-Chowla 1}) and (\ref{Berndt-Chowla 2}) can also be inferred by an argument
of B.C.~Berndt and S.~Chowla (p.~8 of \cite{BerChow}) in combination of their main Theorem therein, applied
with $q=2$.

Property \ref{property 5} implies another interesting fact, already noted in 1979, namely,
\begin{property} \label{property 6}
If $p\equiv 3\pmod 8$, then the number of even quadratic residues $\bmod\,p$ that are $<p/2$ equals $(p-3)/8$.
If $p\equiv 7\pmod 8$, then the number of even quadratic residues that are $>p/2$ equals $(p+1)/8$.
\end{property}
The fact that Property \ref{property 6} is implied by Property \ref{property 5} is noted by 
Emma Lehmer \cite{Lehmer}. 
\\
Finally, we remark that our arguments that led to Property \ref{property 5} furnish another 
expression for $T(p)$, namely,
\begin{equation} \label{T5}
 T(p)=\begin{cases}
  p\left(\jac{(p+1)/4}{p}+\cdots+\jac{(p-3)/2}{p}+\jac{(p-1)/2}{p}\right) & \mbox{if $p\equiv 3\pmod 8$} \\[2mm]
 -p\left(\jac{1}{p}+\jac{2}{p}\cdots+\jac{(p-3)/4}{p}\right) & \mbox{if $p\equiv 7\pmod 8$}
      \end{cases}
\end{equation}


\begin{thebibliography}{9}
%
\bibitem{Apo}\textsc{T.M.~Apostol},\ {\em Introduction to Analytic Number Theory},
  Springer-Verlag, New York 1976.
%
\bibitem{BaCho}\textsc{R.P.~Bambah, S.~Chowla},\ On the sign of the Gaussian sum,
{\em Proc.~Nat.~Inst.~Sci.~India} {\bf 13} (1947), 175-176.
%
\bibitem{Ber}\textsc{B.C.~Berndt},\ Classical theorems on quadratic residues,
{\em Enseign.~Math.} {\bf 22} (1976), 261-304.
%
\bibitem{BerChow}\textsc{B.C.~Berndt, S.~Chowla},\ Zero sums of the Legendre symbol,
{\em Nordisk Mat.~Tidskrift} {\bf 22} (1974), 5-8.
%
\bibitem{BerZah}\textsc{B.C.~Berndt, A.~Zaharescu},\ Finite trigonometric sums and class numbers,
{\em Math.~Ann.} {\bf 330} (2004), 551-575.
%
\bibitem{Chung} \textsc{Kai-Lai Chung},\ Note on a theorem on quadratic residues,
{\em Bull.~Am.~Math.~Soc.} {\bf 47} (1941), 514-516.
%
\bibitem{Cohen} \textsc{H.~Cohen},\ {\em A Course in Computational Algebraic Number Theory},
Springer Graduate Texts in Mathematics No 138, Berlin Heidelberg 1993.
%
\bibitem{Dav}\textsc{H.~Davenport},\ {\em Multiplicative Number Theory - Second Edition},
Graduate Texts in Mathematics Vol.~74, Springer-Verlag, New York 1980.
%
\bibitem{Dick}\textsc{L.E.~Dickson},\ {\em History of the Theory of Numbers}, Vol. III, 
Chelsea Publishing Co., New York 1971.
%
\bibitem{Dir1}\textsc{G.L.~Dirichlet},\  Recherches sur diverses applications de l' analyse infinit\'{e}simale
\`{a} la th\'{e}orie des nombres, {\em Reine Angew.~Math.} {\bf 19} (1839), 324-369.
%
\bibitem{Dir2}\textsc{G.L.~Dirichlet},\ Recherches sur diverses applications de l' analyse infinit\'esimale
\`a la th\'eorie des nombres, seconde partie, {\em J.~Reine Angew.~Math.} {\bf 21} (1840), 134-155.
%
\bibitem{IrRo}\textsc{K.~Ireland, M.~Rosen},\ {\em A classical introduction to modern Number Theory},
Graduate Texts in Mathematics Vol.~84, Springer-Verlag, New York 1982.
%
\bibitem{Hua}\textsc{Hua Loo Keng},\ {\em Introduction to Number Theory}, Springer Verlag, Berlin 1982.
%
\bibitem{Leb}\textsc{V.A.~Lebesgue},\ Suite du Memoire sur les applications du symbole $\jac{a}{b}$,
{\em J.~de Math.} {\bf 15} (1850), 215-237.
%
\bibitem{Lehmer}\textsc{E.~Lehmer},\ Solution of Problem 6156, {\em Am.~Math.~Monthly} {\bf 86} No 2
(1979), 134-135.
%
\bibitem{Lerch}\textsc{M.~Lerch},\ Essais sur le calcul de nombre des classes de formes quadratiques
binaires aux coefficients entiers, {\em Acta Math.} {\bf 29} (1905), 333-424; 
{\em Acta Math.} {\bf 30} (1906), 203-293.
%
\bibitem{Montg}\textsc{H.L.~Montgomery},\ An exponential polynomial formed with the Legendre symbol,
{\em Acta Arith.} {\bf 37} (1980), 375-380.
%
\bibitem{Mo}\textsc{L.~Moser},\ A theorem on quadratic residues, {\em Proc.~Am.~Math.~Soc.} {\bf 2} No 3
(1951), 503-504.
%
\bibitem {Whi}\textsc{A.L.~Whiteman},\ Theorems on quadratic residues,  
{\em Mathematics Magazine} {\bf 23} No 2 (1949), 71-74.
%

%
\end{thebibliography}
\end{document}